
\input amstex

\overfullrule=0pt
\define\ssk{\smallskip}

\define\del{$\partial$}

\documentstyle{amsppt}

\topmatter

\title\nofrills Free genus one knots with large volume \endtitle
\author  Mark Brittenham \endauthor

\leftheadtext\nofrills{Mark Brittenham}
\rightheadtext\nofrills{Free genus one knots with large volume}

\affil   University of  North Texas \endaffil
\address   Department of Mathematics, University of North Texas, Denton, TX 76203  \endaddress
\email   britten\@unt.edu \endemail
\thanks   Research supported in part by NSF grant \# DMS$-$9704811 \endthanks

\abstract In this paper we construct families of hyperbolic 
knots with free genus one, whose complements have arbitrarily
large volume. This implies that these knots have free genus one but 
arbitrarily large canonical genus. \endabstract

\endtopmatter

\heading{\S 0 \\ Introduction}\endheading

A Seifert surface for a knot $K$ in the 3-sphere is an embedded
 orientable surface $\Sigma$, whose boundary equals the knot $K$.
In 1934, Seifert [Se] gave a very simple algorithm for 
constructing a Seifert surface for a knot, from a 
diagram, or projection, $D$ of the knot. Thus every knot 
has a Seifert surface.

Seifert's algorithm always builds a surface whose complement
is a handlebody, something which is known as a {\it free}
Seifert surface. The the minimal genus among all free Seifert 
surfaces for $K$ is called the free genus $g_f(K)$ of $K$, 
while the minimal genus of a surface built by Seifert's 
algorithm applied to a projection of the knot $K$ is 
called the canonical genus $g_c(K)$ of $K$. (In keeping with
this terminology, we will call a surface built by Seifert's 
algorithm {\it canonical}.) The above 
considerations immediately imply that, for any knot, 
$g_f(K)\leq g_c(K)$. It was shown by Kobayashi and Kobayashi 
[KK] that these numbers can be distinct; for $K$ the connected 
sum of $n$ copies of the double of a 
trefoil knot, $g_f(K)=2n$ and $g_c(K)=3n$. 

An unusual feature of these examples is that the free genus 
minimizing surfaces are all compressible. We were interested
in finding examples where the free and canonical genera differ,
but the free genus minimizing surfaces were incompressible. In doing 
so, we were led in a natural way to consider hyperbolic knots, in 
order to exploit a relationship between canonical genus and volume.

In a previous paper [Br] we showed that hyperbolic knots with
bounded canonical genus have complements with bounded volume.
The bound on volume can in fact be chose to be linear in the
canonical genus. In this paper we show, by contrast, that there
is no corresponding bound in terms of the free genus of the knot.

\proclaim{Theorem}There exist hyperbolic knots with free genus 
one and arbitrarily large hyperbolic volume. 
\endproclaim

These two results together show that
there are free genus one hyperbolic knots with arbitrarily
large canonical genus. Since a free genus one Seifert surface for 
a non-trivial
knot is always incompressible, the knots we build also provide examples
of knots with incompressible free Seifert surfaces which cannot
be built by applying Seifert's algorithm to a projection of the 
knot. We should also note that a knot with free genus one has genus one, 
and so we immediately obtain the following corollary.

\proclaim{Corollary} There exist hyperbolic knots with genus one and 
arbitrarily large hyperbolic volume.
\endproclaim


\heading{\S 1 \\ Building free genus one knots}\endheading

It is easy to build knots with free 
genus (at most) one,
simply by building a genus one surface whose spine is an 
unknotted graph in the 3-sphere. The complement of the surface 
is homeomorphic 
to the complement of the graph, and so is a handlebody. 
Consequently, the boundary of the surface has a genus 1
free Seifert surface, and so has free genus at most one.
The best examples of this (and the starting point for our examples)
are the 2-bridge knots (Figure 1) corresponding to the continued fractions
[2u,2v] (see [HT] for notation). These surfaces are in fact isotopic to canonical 
surfaces for different projections of these knots.

This gives an infinite family of free genus one knots. However, 
since all of these knots can be built by doing $1/u$ (and $1/v$) 
Dehn surgeries on two of the unknotted components of a single link
(Figure 1), these knots have hyperbolic volume smaller than the 
hyperbolic volume of the link [Th1], and so have bounded volume.

\ssk

\input epsf.tex

\leavevmode

\epsfxsize=4.9in
\centerline{{\epsfbox{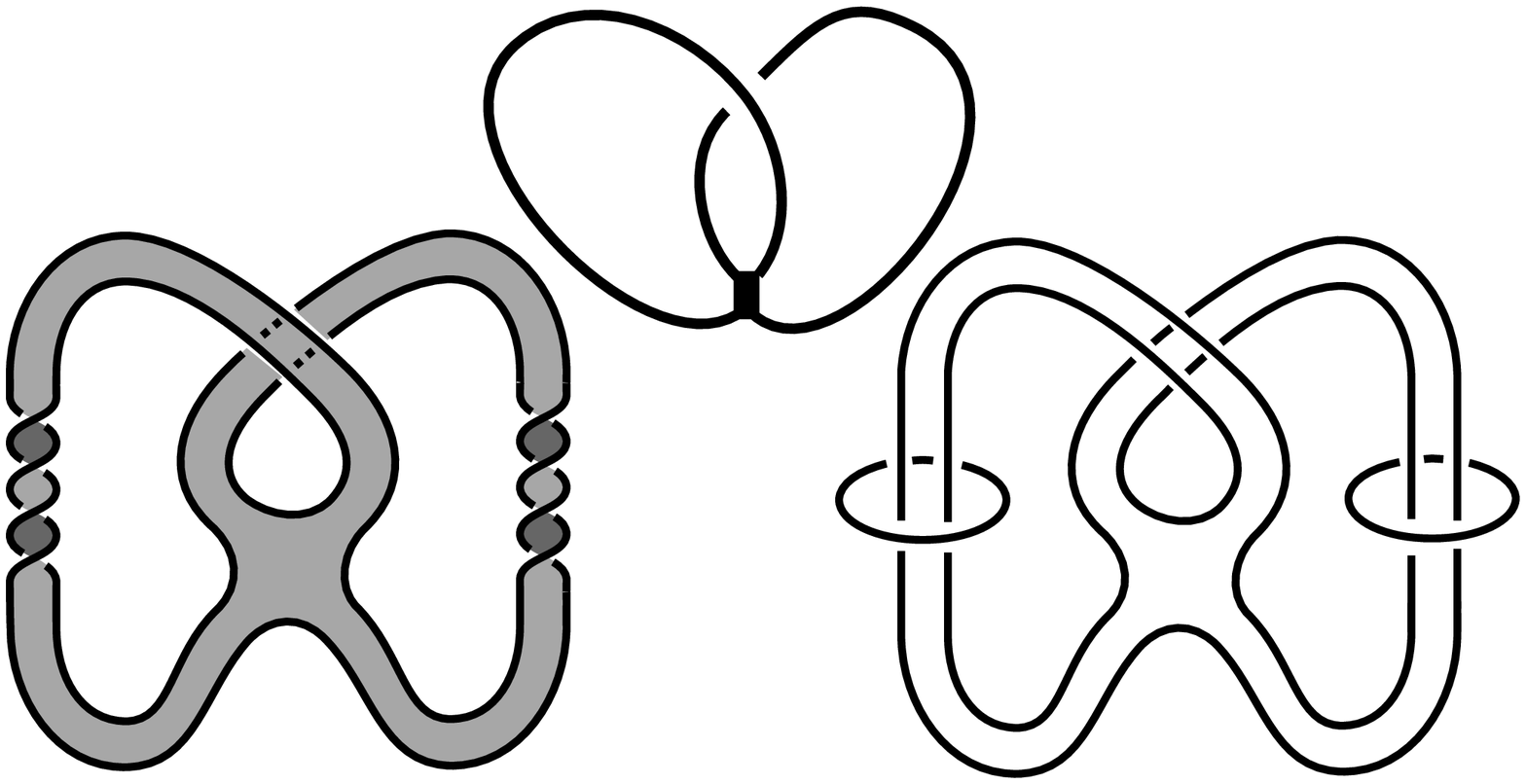}}}

\centerline{Figure 1}

\ssk

In order to insure that our free genus one knots will have 
large volume, we will rely on a result of Adams [Ad], which states 
that a complete hyperbolic manifold with $r$ cusps must have 
volume at least $rV_0$, where $V_0$ is the volume of a regular
ideal tetrahedron. We will therefore build our knots by doing
$1/n_i$ Dehn surgeries on the (unknotted) components of a hyperbolic
link with $r$ components. By a result of 
Thurston [Th1], for 
large values of the $n_i$, the resulting knots will also be
hyperbolic and have volume close to that of the link. What we shall see 
is that for the particular link we choose, the resulting knots can
also be seen to have free genus one Seifert surfaces, and so have 
free genus (at most) one.

The basic idea is to take one of the knots $K$ in Figure 1 and throw an extra 
loop $K_1$ around the 
`waist' of the Seifert surface $F$; see Figure 2. This loop $K_1$ lies on a 2-sphere
$S$ bounding a 3-ball $B$; $B\cap F$ is a disk and $S\cap F$ consists of four
arcs. If we look at $M_0=B\setminus$int($N(F)$), it is a genus four 
handlebody which has the structure of a sutured manifold, where the sutures
are the four loops $S\cap\partial(N(F))$. There are four obvious product disks
for this sutured manifold (sitting in the plane of the paper, in the figure), 
so that, as a sutured manifold, $M_0$ is a product sutured manifold
(four-punctured sphere)$\times I$. In particular, the
four-punctured sphere $B\cap\partial N(F)$ is isotopic,
in $M_0$, to $S\setminus$int$(N(F))$.

\ssk

\leavevmode

\epsfxsize=4in
\centerline{{\epsfbox{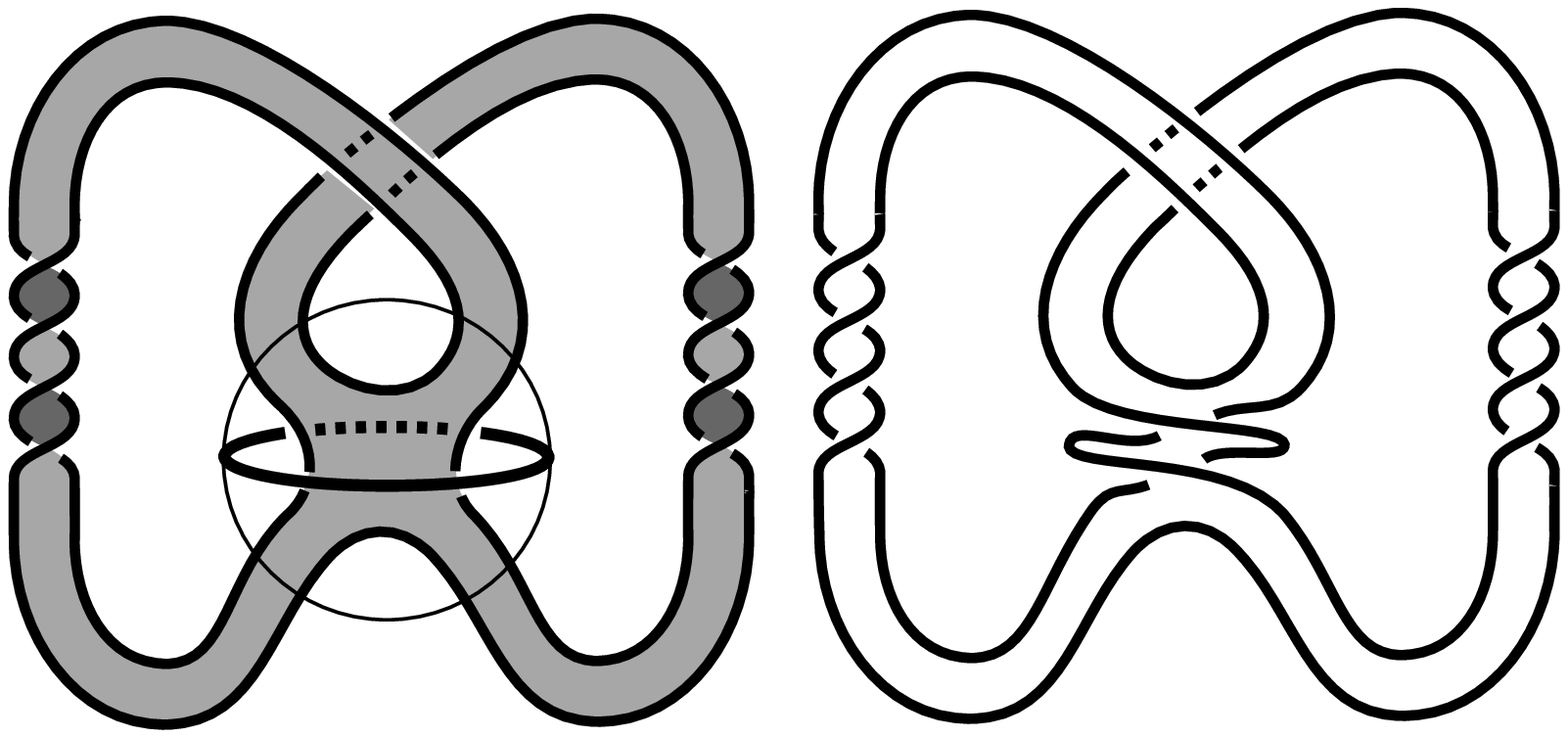}}}

\centerline{Figure 2}

\ssk

If we now do $1/n$ Dehn surgery on $K_1$, then since $K_1$ is unknotted, the 
resulting manifold will
be a 3-sphere, and so the image of $K$ will be a knot $K^\prime$ in the
3-sphere. By Kirby calculus (see [Ro]), $K^\prime$ can be obtained from $K$
by cutting $K$ along a disk $D$ spanning $K_1$, giving the resulting strands $n$ 
right-handed
twists, and regluing. Similarly, the Seifert surface for $K$ gives a Seifert
surface $F^\prime$ for $K^\prime$, by cutting, twisting, and regluing.

However, from the point of view of the ball $B$ (or, more precisely, a ball slightly
smaller than $B$), nothing has really happened; we have cut the ball open along a 
disk, given one half of it $n$ full twists, and reglued {\it by the identity map}. 
Consequently, $B\setminus$int$(N(F^\prime))\cong B\setminus$int$(N(F))$ is also a 
product sutured manifold, and so $B\cap\partial N(F^\prime)$ is isotopic,
in $B\setminus$int($N(F^\prime)$), to $S\setminus$int$(N(F^\prime))=S\setminus$int$(N(F))$.
Therefore, $S^3\setminus$int($N(F^\prime)$) $\cong$ $S^3\setminus$(int($B\cup N(F^\prime)$)
$=$ $S^3\setminus$(int($B\cup N(F)$) $\cong$ $S^3\setminus$int($N(F)$) is also
a handlebody. $F^\prime$ is therefore a free genus one Seifert surface for $K^\prime$.
Note that we could in fact have chosen any loop on $S$ to base this construction on; 
it would still bound a disk in the 3-ball $B$, and so the argument above would go 
through without change. A more complicated loop would, however, unnecessarily
complicate our arguments below.

More generally, we can throw many extra loops $K_i$ around $F$, on concentric 2-spheres, 
alternating which direction around the waist of $F$ we go (Figure 3), and do $1/n_i$ Dehn 
surgeries on each of them. Since without $K$ these loops would together
form a trivial link - they lie on disjoint 2-spheres - the resulting manifold will 
again be the 3-sphere, and so 
$K$ will be taken to a new knot $K^\prime$ in the 3-sphere. By working inductively 
out from the 
centermost 2-sphere, cutting along disks, twisting, and regluing, we can also 
see that our Seifert surface $F$ will be taken to a free
Seifert surface for $K^\prime$; at each step the argument is identical to the one
given above. We therefore can produce knots with free genus at most one by this
iterative construction, as well.

What we do {\it not} yet know is that these knots $K^\prime$ are hyperbolic, or 
that they have
large volume. What we will show, however, is that the links
$L_n$ = $K\cup K_1\cup\ldots K_n$ are all hyperbolic. By the result of Adams, 
for large $n$ they 
therefore have large volume, and so when the $n_i$ are all large, $K^\prime$ will have
large volume (and, in particular, is therefore non-trivial, hence has
free genus exactly one). This will complete the proof of our theorem.

\ssk

\leavevmode

\epsfxsize=3.5in
\centerline{{\epsfbox{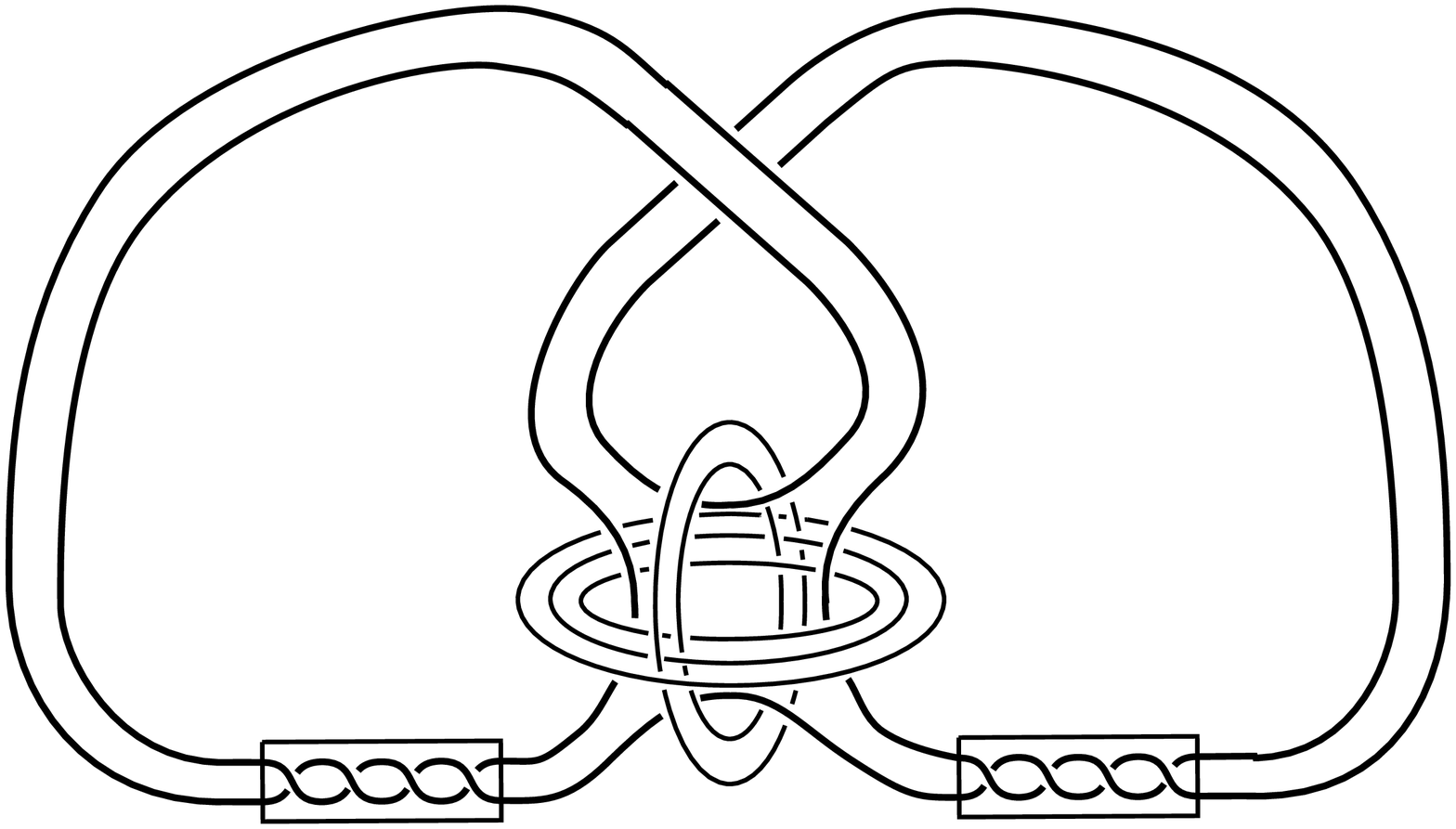}}}

\centerline{Figure 3}

\ssk


\heading{\S 2 \\ The links $L_n$ are hyperbolic: preliminaries}\endheading

We now demonstrate that the links $L_n$ have hyperbolic
complement; that is, the compact manifolds $M_n$=$S^3\setminus$int$(N(L_n)$ 
have hyperbolic interior. By Thurston's Geometrization Theorem [Th2], we 
must show that

\smallskip

(1) $X(L_n)$ = $S^3\setminus$int$N(L_n)$ is irreducible,

\smallskip

(2) $X(L_n)$ is $\partial$-irreducible, i.e., $\partial X(L_n)$ is 
incompressible in $X(L_n)$,

\smallskip

(3) $X(L_n)$ is atoroidal, i.e., an incompressible torus $T$ in 
$X(L_n)$ is parallel to $\partial X(L_n)$, and

\smallskip

(4) $X(L_n)$ is anannular, i.e, any properly embedded incompressible annulus $A$ 
in $X(L_n)$ is $\partial$-parallel.

\medskip

In this section we will set up a few additional assumptions and prove
some preliminary results, which will
allow us to develop the machinery to prove these assertions. The basic idea
is that, since this is a proof by construction, we can (and will) make whatever
assumptions we feel are necessary to bring a wide array of different tools to
bear on the problem, from standard cut and paste arguments to homological 
intersection numbers to normal forms for words in a free group.

The main assumption we will need to make, in order to prove that
the links $L_n$ are hyperbolic, is that our
underlying 2-bridge knot $K$, given by the continued fraction $[2u,2v]$, lying
at the center of the links $L_n$ has $u,v\geq 2$. Experimental
evidence (finding hyperbolic structures using SnapPea [We]) suggests that in 
fact the links are always hyperbolic, without this added assumption, 
but some of our arguments will not go through in greater generality. 
We will also assume that $n\geq 3$, since this will,
in the end, make the verification of condition (4) almost immediate. Again, experimental 
evidence
suggests that this is not a necessary assumption.

Central to our proof that the knots obtained by $1/n_i$ surgeries have free genus at
most one was that the Seifert surface $F$ for the knot $K$ is disjoint from all of the 
added
components $K_i$. In fact, each loop $K_i$ bounds a disk $D_i$ in $S^3$ which 
meets $F$ in an arc $\alpha_i$ (Figure 4). We will assume that we have pushed these disks 
slightly off of one another, so that if $i-j$ is even, then $D_i$ and $D_j$ are disjoint. 
When
$i-j$ is odd and $j<i$, then $D_j\cap D_i$ consists of an arc $\beta_{ji}$ contained
in the interior of $D_i$. In particular, $D_j\cap K_i$ is empty, and $D_i\cap K_j$
consists of two points. Also, for each $i$, $D_i\cap(F\cup D_1\cup\ldots\cup D_{i-1})$
is a finite tree, consisting of parallel arcs $\beta_{ji}$ each pierced by the arc $\alpha_i$
exactly once (see Figure 4). The surface $F$ is two-sided; we will arbitrarily assign it
a normal orientation, and call one side of the surface $F_+$ and the other side $F_-$.
(Formally, we should think of this as being the two sides of $\partial N(F)$ in $X(L_n)$,
but we won't really make such a distinction.)

\ssk

\leavevmode

\epsfxsize=3in
\centerline{{\epsfbox{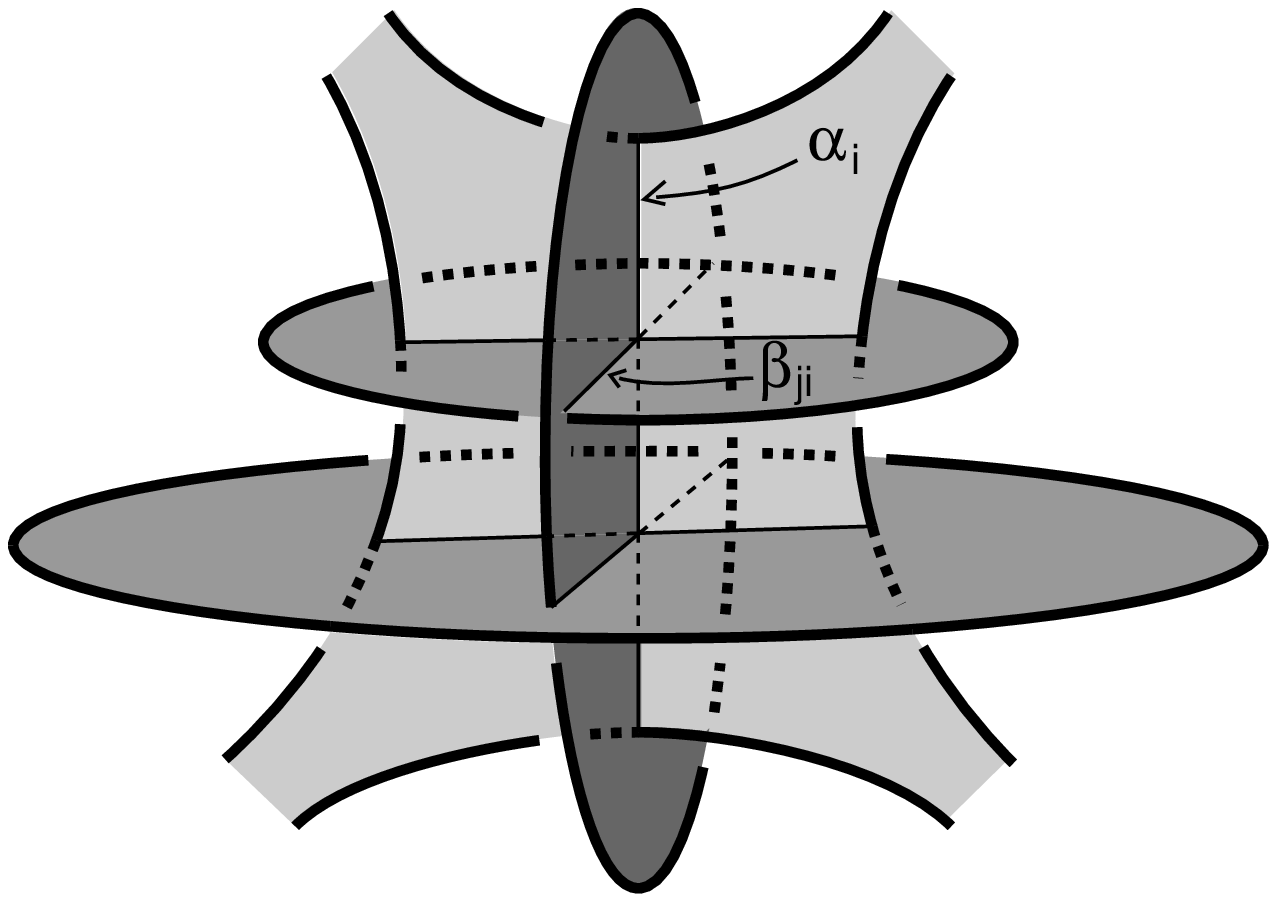}}}

\centerline{Figure 4}

\ssk

An important point to notice is that not only is $S^3\setminus F$ a handlebody (of genus 2),
since $F$ is a free Seifert surface,
but $S^3\setminus (F\cup D_i\cup\ldots\cup D_n)$ is a handlebody (of genus 2), as well.
This is because $A_i$=$D_i\setminus$int$N(F\cup D_1\cup\ldots\cup D_{i-1})$ is an annulus.
Therefore, $X_i$=$S^3\setminus$int$N(F\cup D_1\cup\ldots\cup D_i)$ is
homeomorphic to $X_{i-1}$=$S^3\setminus$int$N(F\cup D_1\cup\ldots\cup D_{i-1})$;
$X_{i-1}\setminus X_i$ is a solid torus neighborhood of $A_i$, so $X_{i-1}$ is obtained from $X_i$
by gluing this solid torus to $\partial X_i$, along an annulus. Pushing $\partial X_i$ to 
$\partial X_{i-1}$ through this solid torus gives an isotopy in $X_{i-1}$ from $X_i$
to $X_{i-1}$.

This fact will allow us to take an inductive approach to our proof of property (3).
We will start with an alleged essential torus $T$, and argue that
we can find a (possibly different) torus disjoint first from $F$, and then, inductively, 
from each 
of the disks $D_i$. After we are done we will have an essential torus disjoint from all of
them, which therefore sits in a handlebody. But since a handlebody is atoroidal, this 
will give us our contradiction.

\ssk

In the rest of this section we collect together several lemmas which will tell us that
certain kinds of intersections of a torus with $F$ and with the disks $D_i$ are not 
possible.

\proclaim{Lemma 1} For every $i$, there is no essential embedded annulus $A$ in 
$S^3\setminus${\rm int}$N(F\cup K_i$) with one \del-component on $F_\pm$ and the
other \del-component on \del$N(K_i)$.
\endproclaim

\demo{Proof} This is easiest to see using a slightly different picture of $F$ (see Figure 5).
$\pi_1(F_+)$, $\pi_1(F_-)$, and $\pi_1(S^3\setminus$int$N(F))=\pi_1(H)$ are all free groups 
of rank two; 
using one of the bases depicted in Figure 5b (depending upon which of 
$F_+,F_-$ we work with), we can see that, as subgroups of $\pi_1(H)=F(a,b)$, 
$\pi_1(F_+)$ is generated by $a^u$ and $ab^v$, and $\pi_1(F_-)$ is generated by 
$b^v$ and $ba^u$. On the other hand, $K_i$ can be represented by either $ab$ or
$ab^{-1}$, depending on the parity of $i$.

\ssk

\leavevmode

\epsfxsize=4.9in
\centerline{{\epsfbox{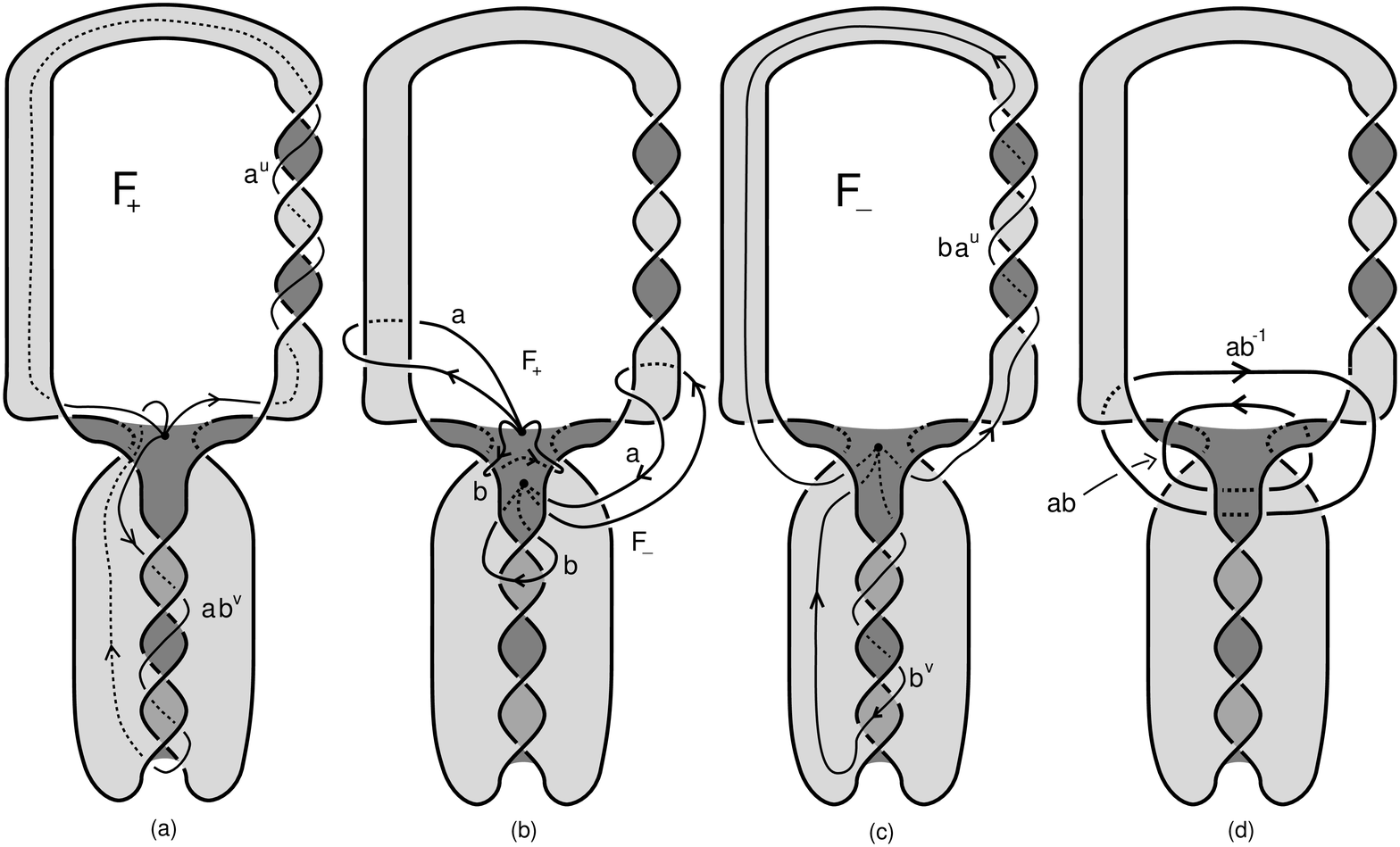}}}

\centerline{Figure 5}

\ssk

From the point of view of homotopy theory, an annulus described by the theorem gives
a free homotopy from a loop representing a power of $ab^{\pm 1}$ to a loop representing 
an element of $\pi_1(F_\pm)$, and 
so from the point of view of fundamental groups, $(ab^{\pm 1})^n$, for some $n$, is
conjugate in $F(a,b)$ to a word in the subgroup
generated by $\{a^u,ab^v\}$ or $\{a^ub^{-1},b^v\}$.

Conjugation preserves exponent sums in a free group, and so in the first case the 
exponent sum for $b$ will be $\pm n=kv$, so $v$ divides $n$, while in the second
case the exponent sum for $a$ will be $n=k^\prime u$, so $u$ divides $n$.
Since by our earlier hypothesis both $u$ and $v$ are greater than 2, we have $|n|\geq 2$,
or $n$=0.
(This is essentially the only place in our proofs where these hypotheses on $u$ and $v$
will be used.)

But $n$=0 implies that $A$ meets \del$N(K_i)$ in a meridian loop (since the
boundary is embedded), and so, capping $A$ off 
with a meridian disk produces a disk $D$ with boundary on $F_\pm$, meeting $K_i$ in a 
single point. Since $F_\pm$ is incompressible, \del$D$ bounds a disk in $F$ which, together 
with $D$ forms a 2-sphere in $S^3$ meeting $K_i$ in a single point, a contradiction (since a 
2-sphere separates $S^3$).
Therefore, $|n|\geq 2$.

\ssk

A word in the free group $F(a,b)$ is said to be in {\it normal form} [MCS] if the 
letters $a$ and $a^{-1}$,
and the letters $b$ and $b^{-1}$, do not occur side by side. Every element of the
free group $F(a,b)$ has a unique normal form, which can be obtained
by starting with a word representing the element and continually cancelling such 
adjacent pairs.

But the word $x(ab^{\pm 1})^nx^{-1}$, when put into normal form, must, since 
$|n|\geq 2$, contain one of the strings

\smallskip

\centerline{$abab$, $baba$, $ab^{-1}ab^{-1}$, $b^{-1}ab^{-1}a$ ,} 

\centerline{$b^{-1}a^{-1}b^{-1}a^{-1}$, 
$a^{-1}b^{-1}a^{-1}b^{-1}$, $ba^{-1}ba^{-1}$, or
$a^{-1}ba^{-1}b$}

\smallskip

\noindent For example, for $x(ab)^nx^{-1}$ with $n\geq 2$, one of the first two strings 
must appear.
This can be proved by induction on the length of the normal word representing $x$.
If we assume $x$ is written in normal form, the only way the initial string $abab$
of the center word, or the final $abab$, can be altered as we shorten our
word to normal form is if $x$ ends in $a^{-1}$, or $x^{-1}$ begins
with $b^{-1}$. But then either $x^{-1}$ begins with $a$ or $x$ ends with $b$, and
so we can write  

\smallskip

\centerline{$x(ab)^nx^{-1}$=$ya^{-1}(ab)^nay^{-1}$=$y(ba)^ny^{-1}$ or}
 
\centerline{$x(ab)^nx^{-1}$=$(z^{-1}b)(ab)^n(b^{-1}z)$=$z^{-1}(ba)^nz$}

\smallskip

\noindent Then by induction (since the word length of the conjugating element has
decreased), we are done; the base case $x$=1 is obvious. The fact that in the inductive step
$ab$ became $ba$ is not a problem, since our conclusion is symmetric in $a$ and $b$;
we simply imagine making our initial statement symmetric in $a$ and $b$ as well.
The other pairs of possibilities 
listed above occur for the other combinations of exponent of $b$ and
sign of $n$. Consequently, the normal form for our word $x(ab^{\pm 1})^nx^{-1}$ contains
both an $a^{\pm 1}$ surrounded on both sides by $b$'s, and a $b^{\pm 1}$ surrounded 
on both sides by $a$'s.

\smallskip

This word is, by our argument above, contained in one of the (free) subgroups 
generated by $\{a^u,ab^v\}$ or $\{ba^u,b^v\}$ . But this is impossible, because $u$ and $v$ are both 
at least 2. In the first case, every occurrence of the letter $b$ will come in the form 
$(b^v)^k$, since $a^u$ contains no $b$'s, and so it is impossible to have a single $b$
surrounded by $a$'s in the normal form for the word. In the second case,
every occurrence of the letter $a$ will come in the form 
$(a^u)^k$, since $b^v$ contains no $a$'s, and so it is impossible to have a single $a$
surrounded by $b$'s in the normal form for the word. Consequently, there can be
no annulus running from \del$N(K_i)$ to $F_\pm$. \qed
\enddemo

\smallskip

\proclaim{Lemma 2} For every $i\neq j$, there is no essential annulus properly embedded in $X(L_n)$ with 
one \del-component on \del$N(K_i)$ and the other \del-component on \del$N(K_j)$
\endproclaim

\demo{Proof} Suppose we have such an annulus $A$; consider $A\cap F\subseteq A$. Since 
$F\cap$\del$A=\emptyset$ and $A\cap$\del$F=\emptyset$, this intersection consists of loops.
Since $F$ is incompressible, any loops which are trivial in $A$ can be removed by disk swapping.
Any remaining loops must all be parallel to \del$A$; if there are any, then an outermost such loop 
cuts off an annulus from $A$ with one \del-component on \del$N(K_i)$ or \del$N(K_j)$ and the other
\del-component on $F_\pm$, contradicting Lemma 1. So $A\cap F=\emptyset$. There are now two 
cases to consider:

\smallskip

\leavevmode

\epsfxsize=3in
\centerline{{\epsfbox{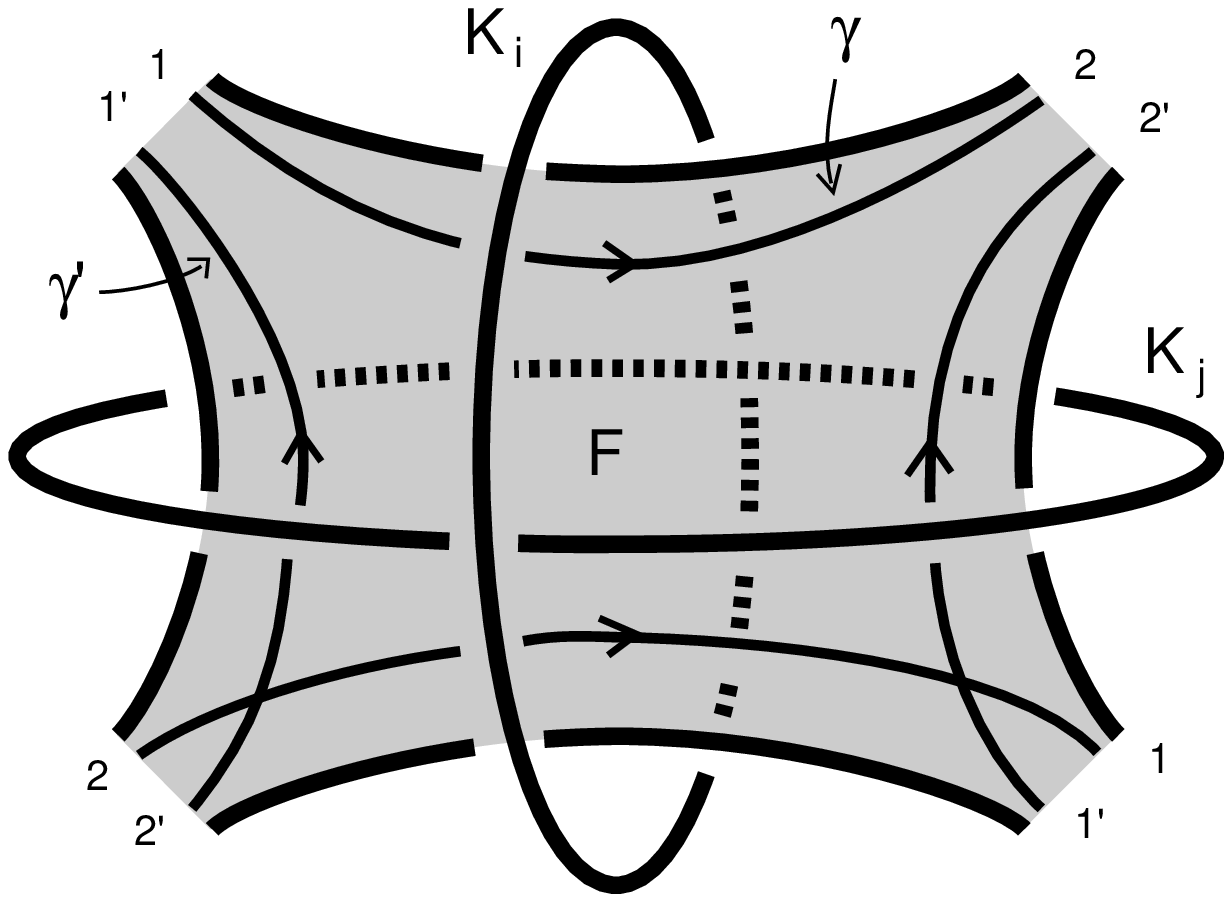}}}

\centerline{Figure 6}

\smallskip

\noindent Case 1: $i-j$ is odd.

\ssk

In this case, $F$, $K_i$, and $K_j$ are as in Figure 6. There is a loop $\gamma$ in $F$ 
which has homological
linking numbers (in $S^3$) 2 with $K_i$ and 0 with $K_j$. \del$A\cap$\del$N(K_i)$ is a
curve of some slope $a_i/b_i$ on \del$N(K_i)$, and \del$A\cap$\del$N(K_j)$ is a
curve of slope $a_j/b_j$ on \del$N(K_j)$. But since $A\cap F=\emptyset$, we have
$A\cap\gamma=\emptyset$, and so $A$ represents a homology in the complement of $\gamma$
between its two boundary curves. But these boundary curves represent the homology 
classes $b_i[K_i]=2b_i$ and $b_j[K_j]=0b_j=0$, and so $b_i=0$. Similarly, a curve 
$\gamma^\prime$ can be found with the appropriate linking numbers, showing that $b_j=0$.
This implies that both \del-components of $A$ are meridian loops; capping off with 
meridian disks gives us a 2-sphere in $S^3$ meeting each of the loops $K_i$,$K_j$ 
exactly once, a contradiction. Therefore, the annulus $A$ cannot exist.

\noindent Case 2: $i-j$ is even.

\smallskip

We may assume $j<i$, and so, setting $k=i-1$, $j<k<i$. 
We then have a situation like in Figure 7. Consider $A\cap D_j\subseteq D_j$. There is an 
arc of $F\cap D_j=\alpha_j$ between the two points of $K\cap D_j$, and since $A$ 
misses $F$, it
misses $\alpha_j$. $A\cap D_j$ therefore consists of trivial arcs, trivial loops, and 
loops surrounding the arc $\alpha_j$ Innermost trivial loops can be removed by disk swapping
with the corresponding disk in $A$, and outermost trivial arcs can be removed since
$A$ is \del-incompressible. Note that this implies that $A\cap$\del$N(K_j)$ misses \del$D_j$, 
and so represents a longitude of \del$N(K_j)$.

\ssk

\leavevmode

\epsfxsize=3in
\centerline{{\epsfbox{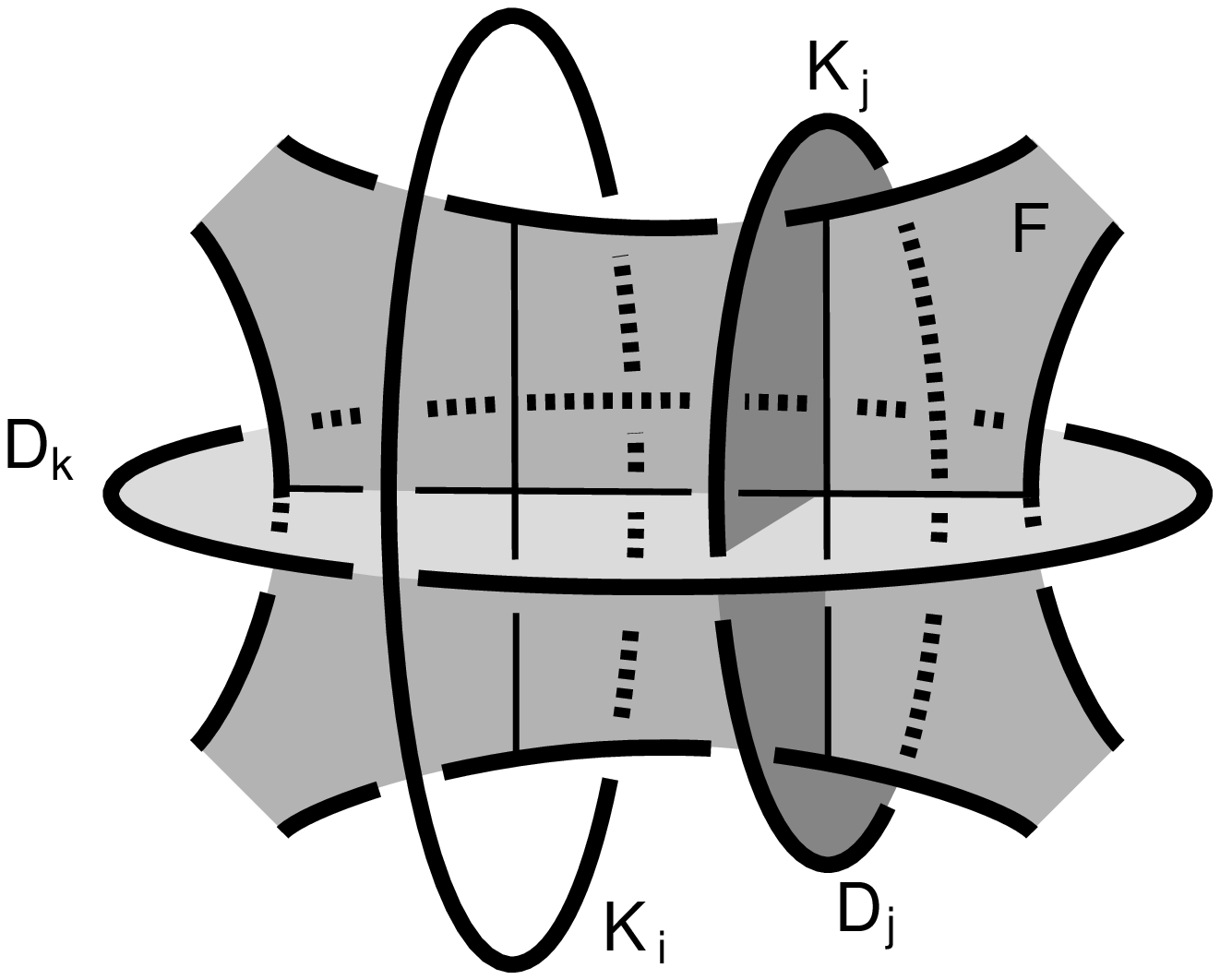}}}

\centerline{Figure 7}

\ssk

This leaves loops travelling around $\alpha_j$; they must 
all be parallel to \del$D_j=K_j$, and therefore all parallel to one another. 
These loops must be non-trivial on $A$, since otherwise we can use the disk bounded by
an innermost trivial loop on $A$, together with the annulus in $D_j$ that the loop
cuts off to build a disk $D^\prime$ disjoint from $F$ with \del$D^\prime=K_j$. But this
contradicts the fact that $K_j$ has homological intersection number 2 with one of the 
two loops $\gamma, \gamma^\prime$ from Figure 6. Then we can use the loop in $A$ 
closest to the \del-component on \del$N(K_i)$, cutting off an annulus from $A$, together
with the annulus in $D_j$ that it cuts off, to build a new annulus $A^\prime$ between
\del$N(K_i)$ and \del$N(K_j)$. Since both boundary components are essential, and lie on
distinct \del-tori, $A^\prime$ is essential. We can then push this annulus off of 
$D_j$ to make it disjoint. Setting $A=A^\prime$, we can therefore assume that
$A\cap D_j=\emptyset$. 

Technically, this new annulus might hit some of the loops $K_r$ for $r<j$; 
what we will actually show, therefore, is that there can be no essential annulus
in $X(K\cup K_j\cup K_k\cup K_i)=X$. This will suffice, since our original annulus $A$
would be essential in this manifold, as well; any compressing disk 
for $A$ in $X$ could be pushed 
off of $F$ by disk swapping, since $A\cap F=\emptyset$, giving us a 
disk $D^\prime$ as in the previous paragraph, a contradiction. Since $A$ 
has its \del-components 
on distinct \del-tori, there can be no \del-compressing disks, either.

\ssk

\leavevmode

\epsfxsize=3in
\centerline{{\epsfbox{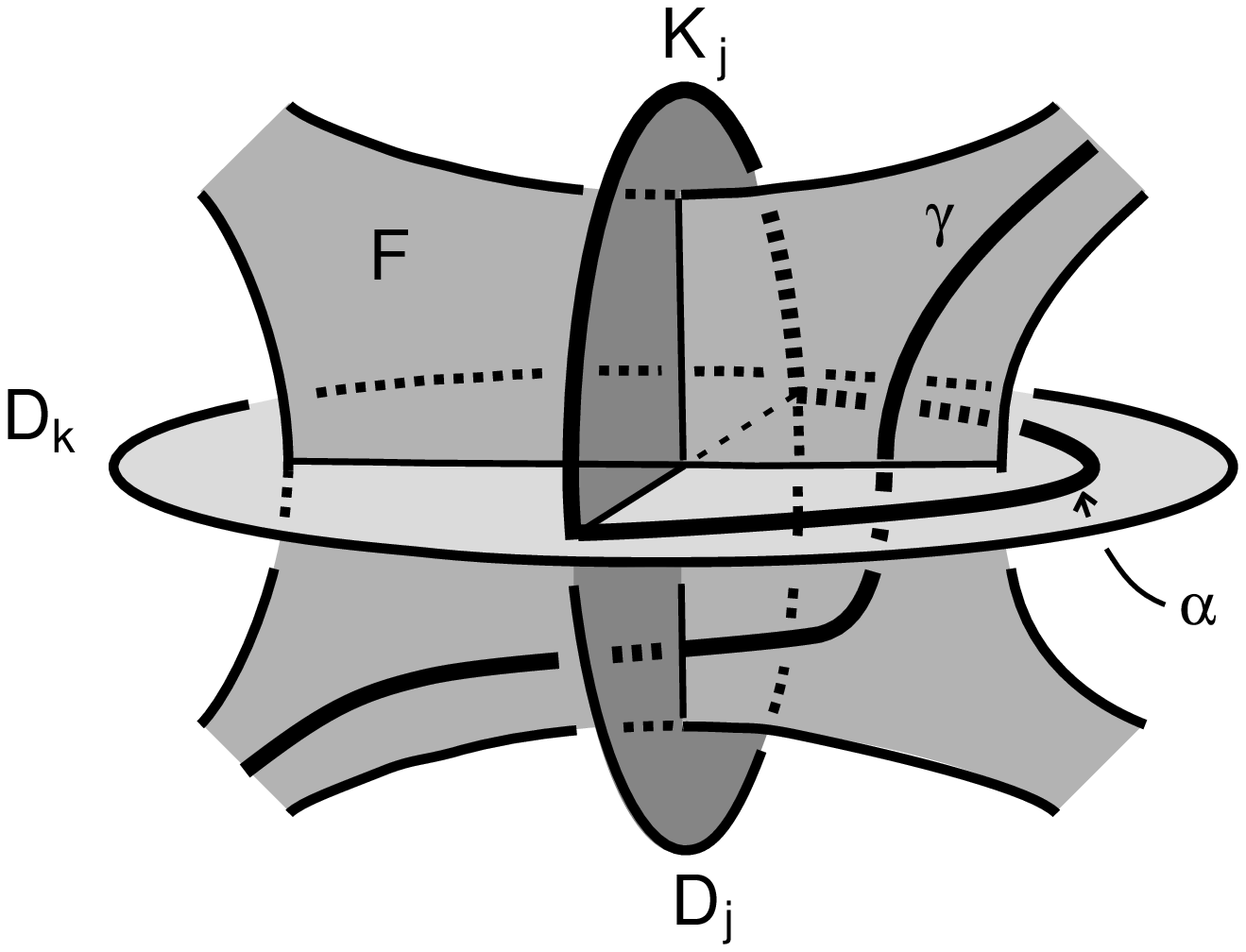}}}

\centerline{Figure 8}

\ssk

Now consider $A\cap D_k$ (Figure 8). $A$ is disjoint from \del$D_k=K_k$, and
since $K_j$ meets $D_k$ in a pair of points, and
$A\cap$\del$N(K_j)$ is a longitude, $A\cap D_k$ consists of circles plus a single arc
$\alpha$  joining the two points of $K_j\cap D_k$. Thinking of 
this arc as being in the annulus $D_k\setminus$int$N(F\cup D_j)$, it is \del-parallel, and so simply goes
to the right or left around the tree $D_j\cap(F\cup D_j)$, say right.

But in $A$, $\alpha$ is \del-parallel, 
since it joins a component of \del$A$ to itself, and so cuts off a disk $\Delta$ from $A$,
which is therefore disjoint from $F$ (since $A$ is). 
$\Delta\cap\partial N(K_j)$ is an arc of a longitude, running above or below the
disk $D_k$, say above. But from the figure, \del$\Delta$ has linking number 1 with a loop
$\gamma$ in $F$; this loop would therefore have to meet the disk $\Delta$, 
a contradiction. \qed
\enddemo


\heading{\S 3 \\ The links $L_n$ are hyperbolic: proofs}\endheading

We now verify the four properties needed to show that the links $L_n$ are
hyperbolic. We work under the assumptions that $n\geq 3$, and the base knot
$K$ has $u,v\geq 2$.

\proclaim{Proposition 1} $X(L_n)$ is irreducible: every embedded 2-sphere bounds a 3-ball.\endproclaim

\demo{Proof} Suppose $S$ is a reducing sphere for $X(L_n)$. $S\subseteq X(L_n)\subseteq S^3$,
and in $S^3$, $S$ bounds a 3-ball $B_1,B_2$ on each side. So we must have $L_n\cap B_i\neq
\emptyset$ for each $i$, otherwise $B_i\subseteq X(L_n)$. One of these 3-balls 
contains $K$, say $B_1$.

$F\subseteq X(L_n)$ is incompressible, since a compressing disk for $F$
in $X(L_n)$ would be a compressing disk in $X(K)$. By a standard 
argument, we can then make $S$ disjoint from $F$: for any innermost loop
of $S\cap F$ in $F$, we can surger $S$ along the corresponding disk in $F$, creating
a pair of 2-spheres, at least one of which must still be a reducing
sphere for $X(L_n)$, with fewer circles of intersection with $F$. Therefore, 
$F\subseteq B_1$, 
since \del$F$=$K\subseteq B_1$. But each component $K_i$ of $L_n$ has non-zero 
linking number with
some loop $\gamma$ on $F$ (Figure 8); in particular, $K_i\cup\gamma$ 
is a nonsplit link. But $K_i\cup\gamma$ is disjoint from $S$, and so is completely
contained in either $B_1$ or $B_2$. Since $\gamma\subseteq F\subseteq B_1$, we
have $K_i\subseteq B_i$ for each $i$. Therefore, $L_n\subseteq B_i$, and 
so $L_n\cap B_2=\emptyset$, a contradiction. So no reducing spheres exist. \qed
\enddemo

\proclaim{Proposition 2} $X(L_n)$ is \del-irreducible: \del$X(L_n)$ is incompressible in $X(L_n)$.\endproclaim

\demo{Proof} Suppose $D$ is a compressing disk for \del$X(L_n)$. Since \del$X(K)$ is incompressible
in $X(K)$ - $K$ is a non-trivial knot - \del$D$ must lie on \del$N(K_i)$ for some $i$. It therefore 
represents a curve of slope $a_i/b_i$ on \del$N(K_i)$. Since
$F$ is incompressible and disjoint from $N(K_i)$, $D$ and $F$ meet in loops
trivial on both, and so we can make $D$ and $F$ disjoint by disk swapping.
But then $D$ is disjoint from the loop $\gamma$ of the previous proof, and so \del$D$ 
is null-homologous
in the complement of $\gamma$. Since \del$D$ represents $b_i[K_i]=b_i$ in $H_1(X(\gamma))$, 
we have $b_i=0$, so
\del$D$ is a meridian loop on \del$N(K_i)$. Capping off with a meridian disk, we get a 2-sphere in 
$S^3$ meeting $K_i$ in a single point, a contradiction. So $D$ does not exist. \qed
\enddemo

\proclaim{Proposition 3} $X(L_n)$ is atoroidal: every incompressible torus in $X(L_n)$ is \del-parallel.
\endproclaim

\demo{Proof} Suppose $T$ is an incompressible torus in $X(L_n)$, and suppose, by way of contradiction,
that it is not \del-parallel. Since $F$ is incompressible in $X(K)$ and disjoint from
$L_n$, it is also incompressible in $X(L_n)$, and by a standard disk swapping argument we can make
$T\cap F$ consist of loops that are essential on both $T$ and $F$. 

Consider $T\cap F\subseteq F$. The loops fall into two types: those that are parallel
to \del$F$, and those that are not (which are all parallel to one another, however, since $F$
is a once-punctured torus). We begin by showing that we can use $T$ to find a
different, essential, torus disjoint from $F$.

Of the \del-parallel loops,
the outermost (i.e., \del$F$-most) loop cobounds an annulus $A$ with \del$F$,
and we can use this annulus to isotope $T$ (in $S^3$; in fact, in $X(L_n\setminus K)$) across \del$F$=$K$ to a torus 
$T^\prime\subseteq X(L_n)$ (see Figure 9). 


\ssk

\leavevmode

\epsfxsize=2.2in
\centerline{{\epsfbox{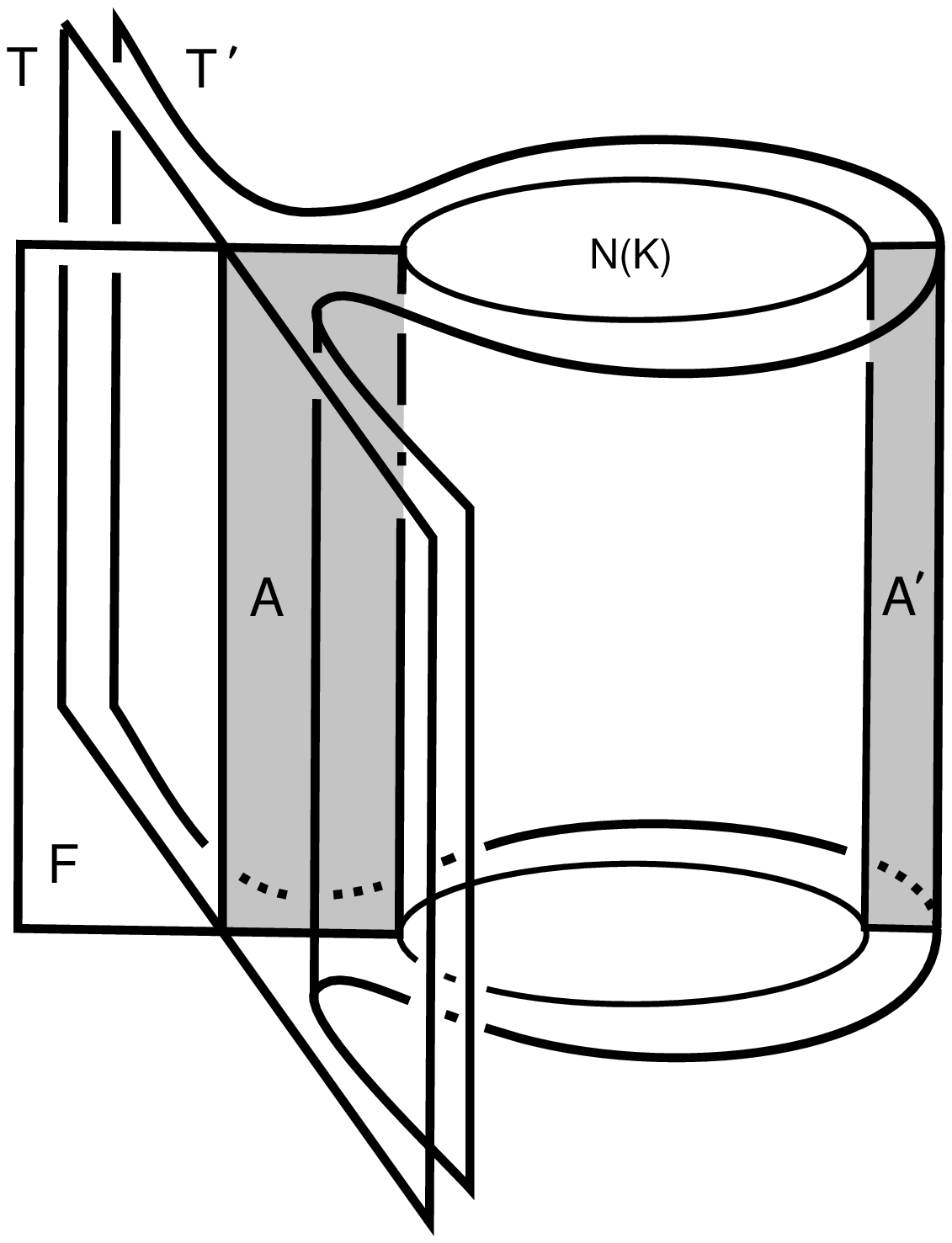}}}

\centerline{Figure 9}

\ssk

\proclaim{Lemma 3} $T^\prime$ is incompressible and not \del-parallel in $X(L_n)$. 
\endproclaim

\demo{Proof} If this is not the case, then one of two things is true:

(1) $T^\prime$ is \del-parallel.

\noindent In this case, if $T^\prime$ is parallel to \del $N(K)$, then
the `dual' annulus joining $T^\prime$ to \del$N(K)$ (Figure 9) would cut the
product region between $T^\prime$ and \del$N(K)$ into a solid torus. By pushing $T^\prime$
back across $K$ using $A^\prime$, we can then see that our original $T$ bounds 
a solid torus, a contradiction. But if $T^\prime$ is parallel to \del$N(K_i)$, then 
$A^\prime$ lies outside of the product region, and the annulus $A^\prime$
together with an annulus in \del$N(K)$ and an annulus in the product region between $T^\prime$ and
\del$N(K_i)$ can be stitched together to form an annulus between \del$N(K)$ and $F$, 
contradicting Lemma 1.
So this case cannot occur.

(2) $T^\prime$ is compressible in $X(L_n)$, via a compressing disk $D$.

\noindent Compressing $T^\prime$ along $D$ produces a 2-sphere $S\subseteq X(L_n)$.
Since $X(L_n)$ is irreducible by Proposition 1, $S$ bounds a 3-ball in $X(L_n)$.
This 3-ball either contains $T^\prime$, or its interior is 
disjoint from $T^\prime$. Therefore, $T^\prime$ either lies in a 3-ball $B_0$ containing $D$,
or bounds a solid torus $M_0$ containing $D$ (Figure 10). 
Also, since $T^\prime\subseteq X(L_n)\subseteq S^3$, $T^\prime$ separates
$X(L_n)$.

The dual annulus $A^\prime$ is incompressible in $X(L_n)|T^\prime$, since \del$A^\prime\cap$\del$N(K)$ is an 
essential loop in \del$N(K)$, hence in $X(L_n)$. If $D$ and $A^\prime$ lie on the
same side of $T^\prime$, then $D\cap A^\prime\subseteq A^\prime$ consists of loops
and arcs. Since $A^\prime$ is incompressible, all of the loops are trivial, and 
since $D$ is disjoint
from \del$N(K)$, no arc joins the two \del-components of $A^\prime$, and so all are 
\del-parallel. By
disk-swapping, we can remove the loops of intersection, and by using the disk cut 
off by an outermost
arc, we can \del-compress $D$ to two disks, at least one of which must still be a 
compressing disk. After
replacing $D$ by one of these disks and continuing, we can eventually find a compressing disk
disjoint from $A^\prime$, which we will still call $D$. 
\ssk

\leavevmode

\epsfxsize=4in
\centerline{{\epsfbox{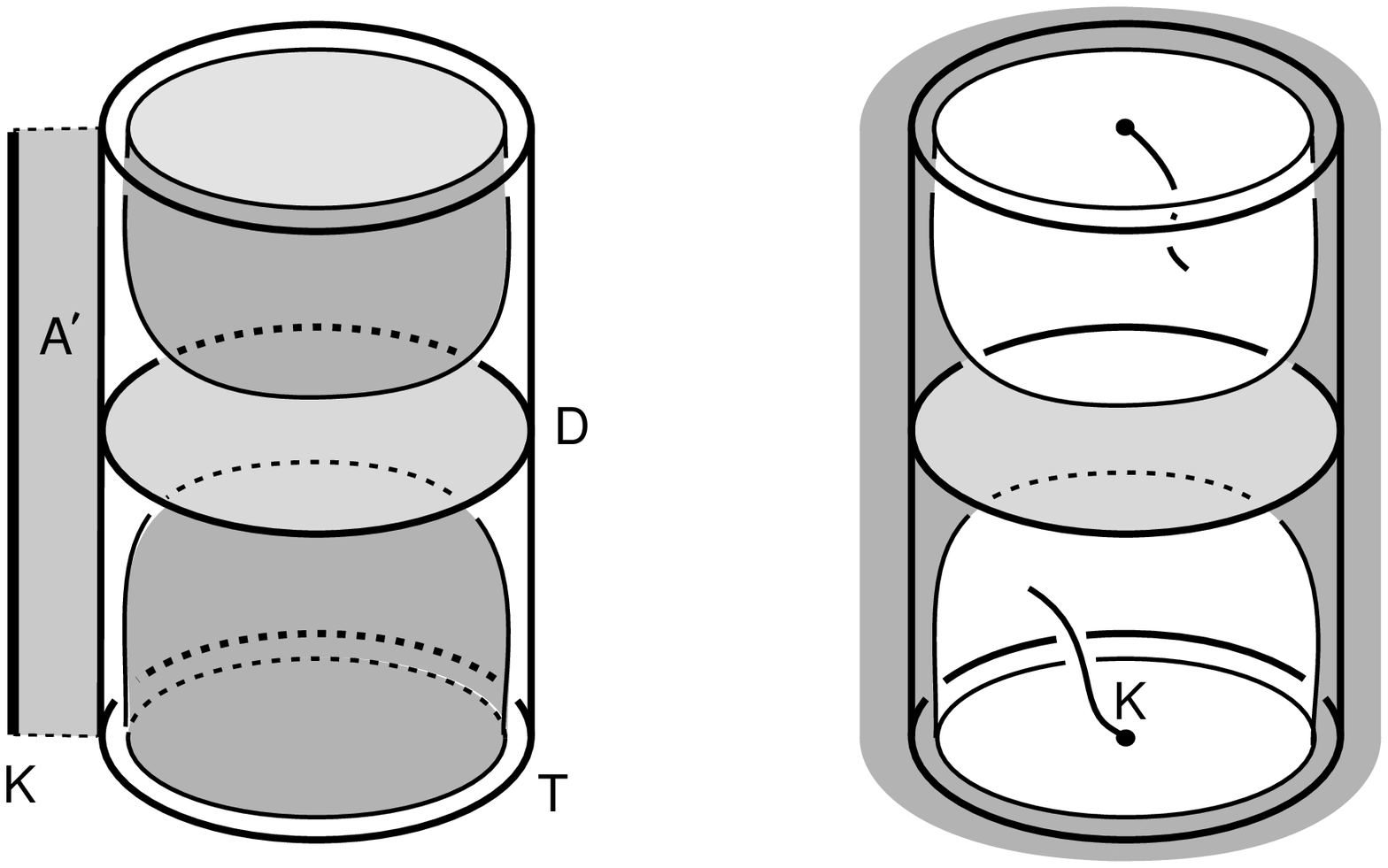}}}

\centerline{Figure 10}

\ssk

If $T^\prime$ bounds $M_0$, then since
$M_0$ must be disjoint from K it is also disjoint from the interior of $A^\prime$ 
(Figure 10), i.e., $A^\prime$ lies outside of $M_0$. The loop
\del$A^\prime\cap T^\prime$ must represent a generator of $\pi_1(M_0)$; otherwise the loop is
meridional and so $A^\prime$ together with a meridian disk of $M_0$ form a compressing disk
for \del$N(K)$ in $X(L_n)$, a contradiction, or the loop represents a non-trivial multiple of 
the generator, and so $A^\prime$ represents an isotopy of $K$ to a non-trivial cable of the core
of the solid torus $M_0$, contradicting the fact [HT] that, for $|u|,|v|\geq 2$, the 2-bridge knots
with continued fraction [$2u,2v$] are hyperbolic.
Therefore, when we push $T^\prime$ back to $T$ 
along $A^\prime$, we see that $T$ is parallel to \del$N(K)$, a contradiction. 

If $T^\prime$ is contained in a 3-ball $B_0$, then since $K$ is disjoint from $B_0$, $D$ and $A^\prime$
lie on the same side of $T^\prime$, and so we may assume,
by the above argument, that they are disjoint. But then when we
push $T^\prime$ back across \del$N(K)$ via $A^\prime$, the compressing disk $D$ persists, so 
$T$ is compressible, a contradiction. So this case also cannot occur. \qed
\enddemo


Therefore, pushing $T$ across $A$ to $T^\prime$ will always result in another incompressible,
non-\del-parallel torus in $X(L_n)$. Continuing for all of the \del-parallel loops
of $T\cap F$, we arrive at a new essential torus, which we will still call $T$, having no
loops of intersection with $F$ parallel to \del$F$, i.e., all loops of intersection
are non-separating on $F$.
These loops cut $T$ into annuli, and $F$ into annuli and a once-punctured annulus. Since
$T$ is separating, none of the annuli in $F$ run from one side of $T$ to the other. One of 
the boundary components of the once-punctured annulus is the longitude of \del$N(K)$.
Note that $F$ and $T$ cannot meet in a single loop $\gamma$; since such a loop is non-separating on $F$,
a loop in $F$ meeting $\gamma$ in a single point is a loop meeting $T$ in a single point,
implying that $T$ is not separating. In fact, this implies that $F$ and $T$ meet in an even
number of loops, since otherwise the same loop will meet $T$ in an odd number of points,
implying they have non-zero homological intersection number, so $T$ is non-trivial in $H_2(S^3)$,
a contradiction. $F|T$ therefore consists of an odd number of annuli and a once-punctured annulus.

Because $X(K)$ is hyperbolic, $T$, thought of as sitting in $X(K)$, must be either compressible or
\del-parallel. It cannot be \del-parallel, however, since then the once-punctured annulus 
component $P$ of $F|T$ would have to live in the product region $T\times I$ between $T$ and \del$N(K)$.
$P$ must therefore be compressible in $T\times I$ (since its fundamental group, being
non-abelian, cannot inject), implying that $F$ is compressible in $X(K)$, a contradiction.

We also know that since $T\subseteq X(L_n)\subseteq X(K)\subseteq S^3$, $T$ bounds a solid 
torus $M_0$ in $S^3$. In fact, $T$ must either bound a solid torus in $X(K)$, or $K$ must
lie in a 3-ball in a solid torus with boundary $T$. 
For if not, then $T$ cannot bound a solid torus on both sides ($K$ would 
be disjoint from one of them).
$T$ is therefore the boundary of a neighborhood of a non-trivial knot, and so is incompressible 
on the side away from $M_0$. Therefore
$K$ lies in $M_0$ and $T$ is incompressible on the side away from $K$. But $K$ is 
not isotopic in $M_0$ to the core $C$ of $M_0$, since $T$ is 
not \del-parallel in $X(K)$, and $C$ is a non-trivial knot (since $T$ is the 
incompressible boundary of $X(C)$). So either $T$ is incompressible in 
$X(K)\cap M_0= M_0\setminus$int$N(K)$, so
$K$ is a satellite knot (and therefore not hyperbolic, a contradiction), or 
$T$ is compressible in $X(K)\cap M_0$, and so $K$ misses a meridian disk for $M_0$ 
(the only compressing disk we could have).
So either $M_0\subseteq X(K)$, or $K\subseteq M_0$ and misses a meridian disk for $M_0$.

If $M_0\subseteq X(K)$ and $T\cap F\neq\emptyset$, then there must be at least one annulus of $F|T$ 
in $M_0$, so there exists
an outermost such annulus $A$. If we split $T$ open along $A$ and glue two parallel copies of $A$ onto 
the resulting annuli, we obtain two new tori $T_1$ and $T_2$ in $X(L_n)$, each bounding a solid sub-torus
in $X(K)$ (see Figure 11), and joined to $F$ by `dual' annuli $A_1$ and $A_2$. Neither of these tori 
can be \del-parallel in $X(L_n)$, otherwise the argument of Lemma 3 will find an annulus contradicting 
Lemma 1. If $T_1$ is compressible, by a compressing disk $D_1$, then since $T_1$ is disjoint from $F$
we can by the usual disk-swapping process make $D_1$ disjoint from $F$ as well. 
If $D_1$ and $A_1$ lie on the same side
of $T_1$, then $D_1\cap A_1\subseteq A_1$
does not contain an essential loop, since otherwise the sub-annulus cut off in $A_1$ parallel to $F$ 
together with the the subdisk cut off in $D_1$ would give a (singular) compressing disk 
for $F$, a contradiction. Therefore by the same process 
used in the proof of Lemma 3, we can make $D_1$ disjoint from $A_1$ as well. 

\ssk

\leavevmode

\epsfxsize=2in
\centerline{{\epsfbox{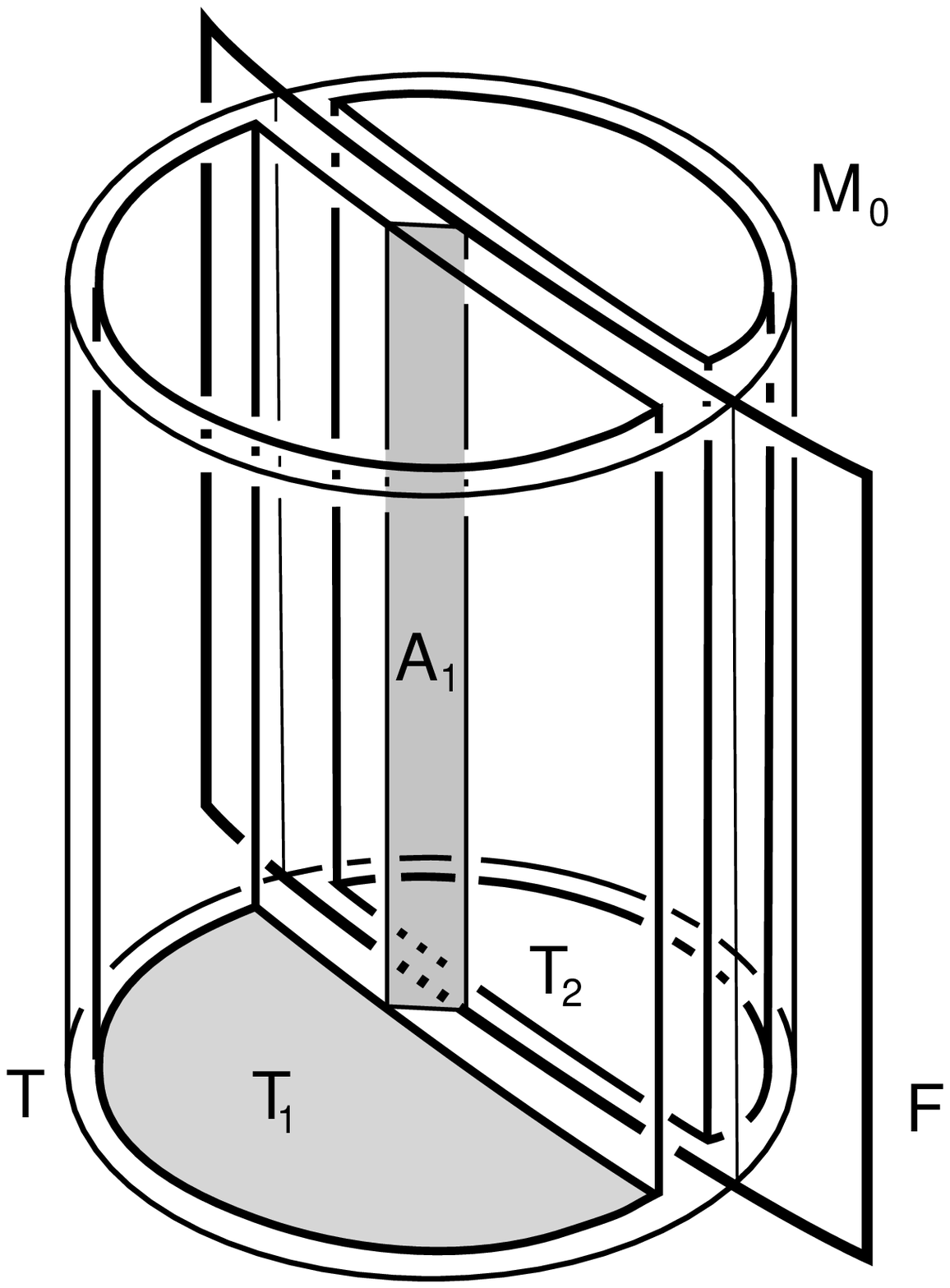}}}

\centerline{Figure 11}

\ssk

But then, as in Lemma 3, either $T_1$ bounds a solid torus $M_1$ in $X(L_n)$, and $A_1$ lies on the 
opposite side of $T_1$ (since $F$ does), or $T_1$ lies in a 3-ball $B_1$, and $D_1$ and $A_1$ lie on 
the same
side of $T_1$. But in the first case, this solid torus must then be identical to the one we see inside
of $M_0$ in Figure 11, and so we can use this solid torus to isotope the annulus of $T|F$ in $T_1$
across $F$, reducing the number of components of $T\cap F$. In the second case, since $D_1$ and $A_1$ 
are disjoint, $D_1$, $A_1$, and an annulus in $T_1$ between their \del-components together
form a compressing disk for $F$, a contradiction. In this case, therefore, we can always
reduce the number of components of $T\cap F$.

\ssk

If $K\subseteq M_0$ and misses a meridian disk $D$ for $M_0$, and $T\cap F\neq\emptyset$, then by 
the incompressibility of $F$
and the irreducibility of $X(K)$ we can isotope $D$ rel \del$D$ to remove any circles of intersection with
$F$. $D\cap F$ must still then be non-empty, since otherwise $D$ together with an annulus in $T$ between
\del$D$ and a component of $T\cap F$ gives a compressing disk (in $X(K)$) for $F$. In particular, 
 the loops of $T\cap F$ are not meridians for $M_0$; otherwise, by isotopy we 
could make \del$D$ (and therefore $D$) disjoint from $F$. By an isotopy of $F$ in $X(L_n)$, 
we can assume that these loops meet \del$D$ minimally.

If $T\cap F$ has only two components, then $F\cap M_0$ consists of a once-punctured annulus $P$, with
two \del-components on \del$M_0$, and the other \del-component equal to $K$. 
Loops of $P\cap D\subseteq D$ can be removed by disk swapping, since $P\subseteq F$ is
incompressible in $X(K)\cap M_0$; $P\cap$\del$M_0$ consists of loops essential in $F$.
Then $P\cap D$ consists of arcs; 
an outermost such arc $\alpha$ must join distinct \del-components of $P$ together, since $T|P$ consists 
of a pair of 
annuli, each of whose \del-components come from distinct components of \del$P$. But then the outermost
disk this arc cuts off, together with the annulus in $T$, can be used to build a disk $D^\prime$
with \del$D^\prime\subseteq F$ (Figure 12). Since $F$ is incompressible, \del$D^\prime$ bounds a disk in $F$.
But since $\alpha$ joins distinct components of \del$P$, $\alpha$ cuts $P$ into an annulus, one
of whose \del-components is \del$D^\prime$, implying that $F$ is a disk, a contradiction. So
$T\cap F$ must have more than two components.

\ssk

\leavevmode

\epsfxsize=2.5in
\centerline{{\epsfbox{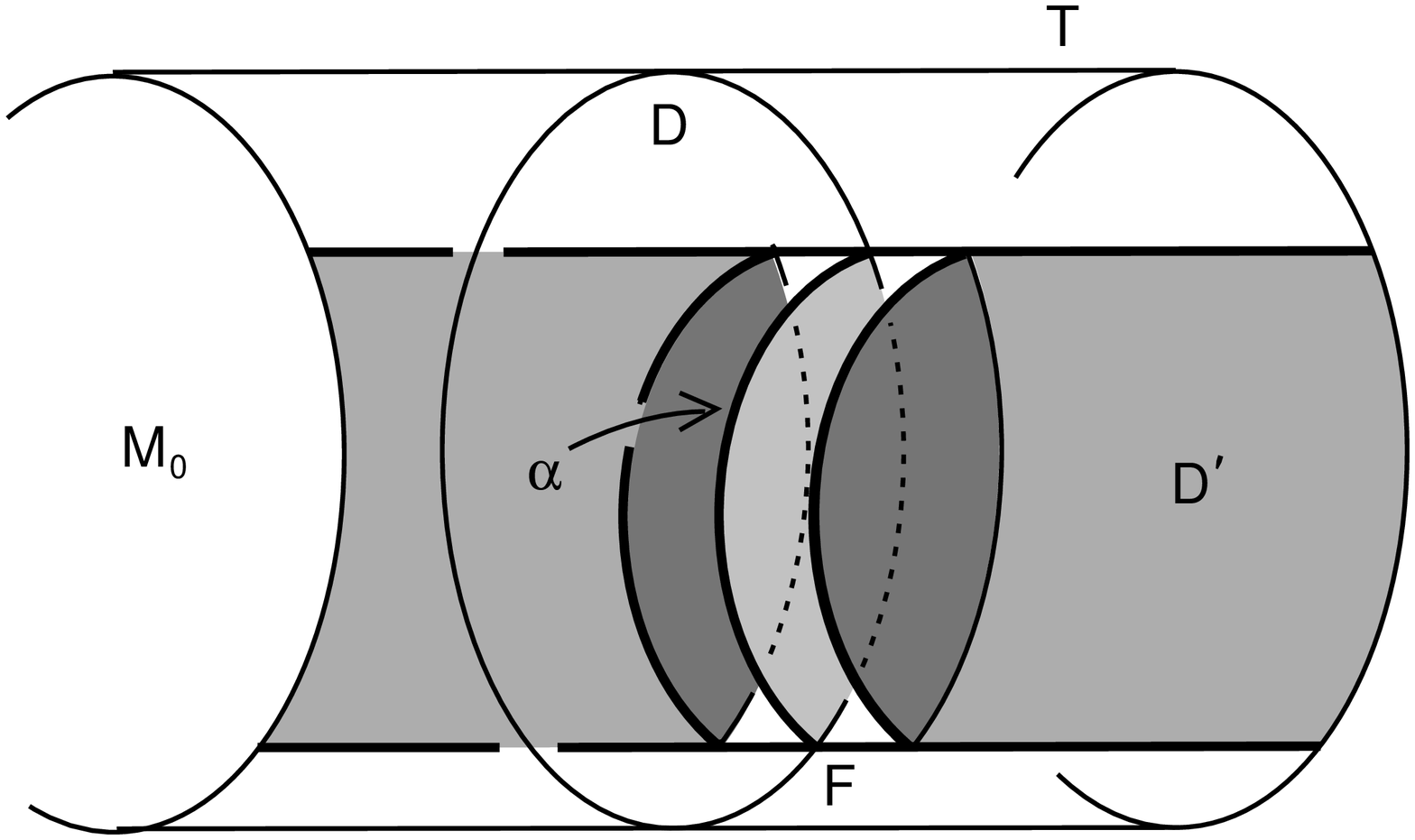}}}

\centerline{Figure 12}

\ssk

If $T\cap F$ consists of more than two components, then there is at least one annulus component of 
$F\cap M_0$. An outermost such annulus cuts a solid torus $M_1$ off from $M_0$. If $M_1$ does not 
contain $K$, then we can apply the argument given above
for the case that $M_0\subseteq X(K)$ to show that we can reduce the number of 
components of $T\cap F$ by an isotopy of $T$. If $M_1$ does contain $K$, then it also contains
$P$, and so we have a situation 
identical to the one in the previous two paragraphs; the meridian disk for $M_1$ is a sub-disk of $D$, 
and so misses $K$. This will lead us to the same contradiction.

\ssk

Therefore, we can either replace $T$ with an essential torus $T_1$ disjoint from 
$F$, or reduce the number of components of $T\cap F$ by an isotopy of $T$ in $X(L_n)$.
Eventually, therefore, we can find an essential torus $T$ with $T\cap F=\emptyset$.

\ssk


Once we have found an essential torus $T$ with $T\cap F=\emptyset$, we turn our 
attention to
$T\cap D_1$, where $D_1$ is the disk of Section 1 bounding $K_1\subseteq L_n$. 
Since $T$ is disjoint from $F$ and $K_1$, after removing 
trivial circles of intersection in $D_1$, which are
therefore trivial on $T$, as well, $T\cap D_1$ consists of loops which 
miss the arc $F\cap D_1$, so all of the loops are parallel, in $D_1$, to \del$D_1$=$K_1$.
The outermost such loop cuts off an annulus $A_1$ from $D_1$ (Figure 13), which we use 
as in Lemma 3 to push $T$ across $K_1$ to a new torus $T^\prime$ in $X(L_n)$.

\ssk

\leavevmode

\epsfxsize=2.5in
\centerline{{\epsfbox{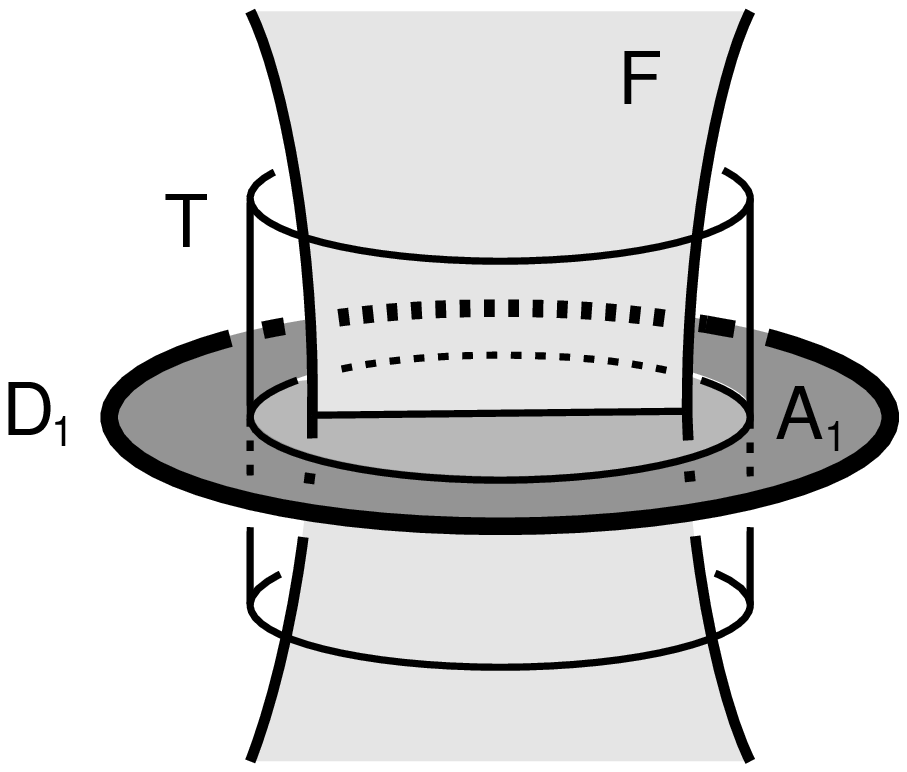}}}

\centerline{Figure 13}

\ssk

\proclaim{Lemma 4} $T^\prime$ is incompressible and not \del-parallel in $X(L_n)$.
\endproclaim

\demo{Proof} Most of our arguments follow the same line as the proof of Lemma 3.
If $T^\prime$ is parallel to \del$N(K_i)$, then if $i\neq 1$, we can use the dual
annulus $A_1^\prime$ from $T^\prime$ to \del$K_1$ together with an annulus 
in the product region
to build an annulus in $X(L_n)$ between \del$N(K_1)$ and \del$N(K_i)$, contradicting
Lemma 2. If $T^\prime$ is parallel to \del$N(K)$, then $A_1^\prime$ together with an
annulus in the product region and an annulus in \del$N(K)$ gives an annulus 
in $X(L_n)$ between \del$X(K_1)$ and $F$, contradicting Lemma 1. And if $T^\prime$
is parallel to \del$N(K_1)$, then the dual annulus $A_1^\prime$ splits the product
region into a solid torus; pushing $T^\prime$ back to $T$ along $A_1^\prime$ 
essentially preserves this solid torus, implying that $T$ bounds a solid torus in $X(L_n)$,
a contradiction. So $T^\prime$ is not \del-parallel.

If $T^\prime$ is compressible, then either $T^\prime$ bounds a solid torus $M_1$, or
$T^\prime$ is contained in a 3-ball $B_1$ in $X(L_n)$. If $T^\prime$ bounds $M_1$, then
 as before $M_1$ is disjoint from the interior of $A_1^\prime$. Note also that
$M_1$ is disjoint from $F$, since $T\cap F=\emptyset$ and $M_1$ does not meet \del$F$=$K$.
As before, \del$A_1^\prime\cap M_0=\gamma_1$ must represent a generator of $\pi_1(M_1)$. 
Otherwise, either $\gamma$ is a meridian $M_1$, and a meridian disk together with 
$A_1^\prime$ gives a disk with boundary $K_1$ disjoint from $F$, implying that $K_1$ has
linking number zero with every loop in $F$, a contradiction, or $\gamma_1$ represents
a non-trivial multiple of the core $C$ of $M_1$, so $\gamma_1$, and therefore $K_1$, is 
homologous to a multiple $r$ of $C$ in the complement of $F$, implying that $K_1$ 
has linking number a multiple of $r$ with every loop in $F$, contradicting the fact
that it has linking number one with some loops (e.g., the $\gamma$ of Figure 8 above).
But now when we push $T^\prime$ back to $T$ using $A_1^\prime$, $K_1$ is pushed to the core
of the solid torus $M_1$, implying that $T$ is \del-parallel, a contradiction.

Finally, if $T^\prime$ and its compressing disk $D$ lie in a 3-ball $B_1$ (the only 
remaining possibility), then $D$ and $A_1^\prime$ lie on the same side of $T^\prime$,
and so, as in Lemma 3, we can make $D$ disjoint from $A_1^\prime$. Then when we push 
$T^\prime$ back to $T$ across $A_1^\prime$, the compressing disk persists, so $T$ is
compressible, a contradiction. So $T^\prime$ must be incompressible and not \del-parallel.
\qed
\enddemo

We can apply this argument to each loop of $T\cap D_1$ in turn, pushing them across
$K_1$ to obtain a new essential torus. After carrying this out for all loops
of intersection, we can then assume that $T\cap F$ and $T\cap D_1$ are both empty. We then
turn our attention to $T\cap D_2$, which as before consists of loops in $D_2$ parallel
to \del$D_2$=$K_2$. By the same process as in Lemma 4, we remove these loops of 
intersection as well. Continuing, we eventually find an essential torus $T$ which is disjoint
from $F$ and all of the disks $D_i$, $i=1,\ldots,n$. 
$T$ therefore lives in $S^3\setminus$int$(F\cup D_1\cup \cdots\cup D_n)$, which as we
remarked in Section 2, is a handlebody $H$. But every torus in a handlebody is
compressible ($\pi_1(T)$ is not free, so it cannot inject into $\pi_1(H)$). This 
compressing disk misses $F$ and all of the $D_i$, and so it lives in $X(L_n)$.
Therefore $T$ is compressible in $X(L_n)$, a contradiction. So no such (original)
torus can exist; $X(L_n)$ is atoroidal. \qed

\enddemo


\proclaim{Proposition 4} For $n\geq 3$, $X(L_n)$ is anannular: every incompressible 
annulus is \del-parallel.
\endproclaim

\demo{Proof} The argument here is standard, we simply use the facts that $X(L_n)$
is irreducible and atoroidal, and has at least four \del-components (since $n\geq 3$).

If $A$ is an incompressible annulus, then if $A$ runs between distinct \del-components
$T_1,T_2$, then $T$=\del$N(A\cup T_1\cup T_2)\setminus(T_1\cup T_2)$ is a torus which
separates pairs of \del-components of $X(L_n)$, so cannot be \del-parallel, and 
must therefore be compressible. But a compressing disk will split $T$ into a 2-sphere
which also separates components of $L_n$, implying that $X(L_n)$ is reducible,
a contradiction. So \del$A$ is contained in a single \del-component $T_1$. Then 
$T$=\del$N(A\cup T_1)\setminus T_1$ is a torus which separates $T_1$ from at least 
three other \del-components. So if $T$ is \del-parallel, it is parallel to $T_1$,
so $A$ lives in a product $T\times I$, and so is \del-parallel. If $T$ is compressible, then
the compressing disk splits $T$ into a 2-sphere separating $T_1$ from 
at least three other \del-components, giving a reducing sphere for $X(L_n)$, a
contradiction.

The only possibility which does not lead to a contradiction, therefore, is that $A$ is 
\del-parallel. Therefore, $X(L_n)$ is anannular. \qed

\enddemo

\ssk

\heading{\S 4 \\ Concluding remarks}\endheading

With this we have finished our proof that the links $L_n$ are hyperbolic. By applying the
construction of Section 1, we can therefore build infinitely many (distinct) knots
with free genus one and volume larger than any fixed constant. We find it both 
amusing and embarrassing to note that we can, however, not exhibit a single explicit
example of this phenomenon, for any fixed constant. 
Existence of our examples is guaranteed only by Thurston's
hyperbolic Dehn surgery theorem, which provides no explicit estimate the sizes of
coefficients $n_i$ sufficient to guarantee hyperbolicity of the knots we build.
And while there are estimates of the volume of a hyperbolic manifold after Dehn filling 
(see, e.g., [NZ]), these are asymptotic estimates, giving no explicit lower
bounds in terms of the $n_i$.

\smallskip

\leavevmode

\epsfxsize=4.9in
\centerline{{\epsfbox{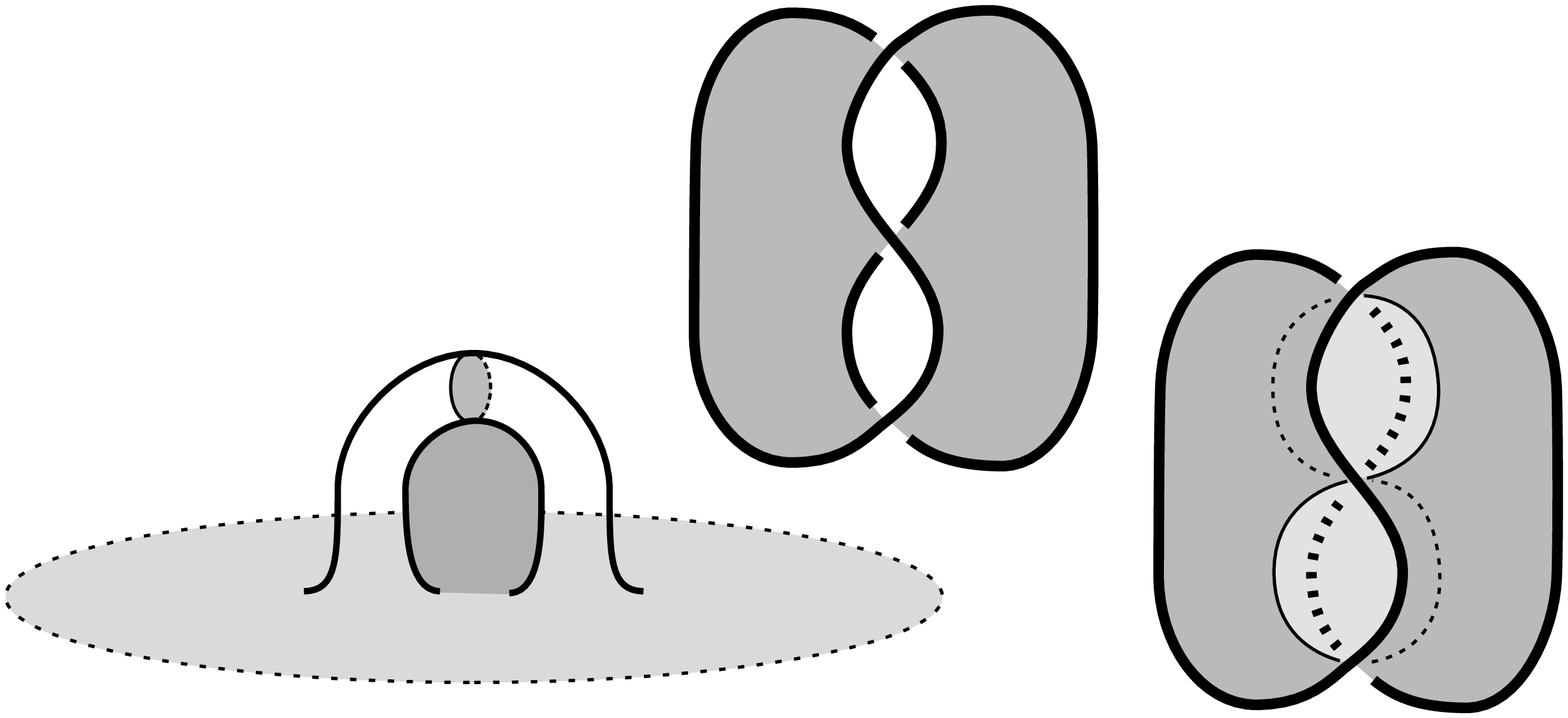}}}

\centerline{Figure 14}

\smallskip

It is clear that similar constructions can be carried out starting with free genus one
knots other than the [$2u,2v$] 2-bridge knots. All that is really needed to carry over our
proofs is an analogue to Lemma 1. Our interest here was in finding hyperbolic knots with
free genus one; we should note that Ozawa [Oz] has, on the other hand, determined all
of the satellite knots with free genus one. In that paper, he conjectures that any 
free genus one knot whose exterior contains an essential closed surface must also 
have tunnel number one. His paper essentially consists of a proof of this conjecture
when the exterior contains an essential torus. This, together with a result of Goda
and Teragaito [GT], gives the required characterization.

A free Seifert surface remains free after stabilization, i.e., 
after adding a trivial 1-handle the surface (see Figure 14).
It is easy to see that a canonical 
surface stabilizes to a canonical one; stabilization can be
thought of as boundary connected sum with a canonical surface
for the unknot (Figure 14), and the connected sum of a diagram for
$K$ and a diagram for the unknot is a diagram for $K$. 
It would be interesting to determine whether or not 
any two free Seifert surfaces are stably equivalent, i.e., they 
become isotopic after a sufficient number of stabilizations. (Since 
the effect, on the complement, of a stabilization is to boundary-connect 
sum with two solid tori, you need to start with handlebody 
complement in order to get handlebody complement.) This is 
probably not unreasonable, since this operation is very similar 
to the stabilization of Heegaard decompositions (see, e.g., [AM]), where stable
equivalence is known.

This stabilization process raises several further interesting 
possibilities. Can every free Seifert surface be stabilized to a 
canonical one? In particular, might it be the case that
any free Seifert surface, if it is stabilized to have genus
equal to (or greater than) the canonical genus, must 
then always be canonical?
This would imply that if the free genus equals the canonical genus, then
every free Seifert surface is canonical.
The standard conjecture in the theory of Heegaard decompositions
(see, e.g., [Ki, Problem 3.89]) seems to be that for any pair of 
Heegaard decompositions of a 3-manifold, they are equivalent
after stabilizing to a genus one higher than the larger of the two.
Perhaps a similar result is true for free Seifert surfaces, as well.

\smallskip

\centerline{\bf Acknowledgements}

\smallskip

The author wishes to thank Jeff Weeks for many helpful conversations during the
course of this work, and Tsuyoshi Kobayashi and Nara Women's University for 
providing the setting in which these conversations took place. The author also 
wishes to thank Susan Hermiller for help in deciphering normal forms in free groups.


\medskip
\Refs

\refstyle{A}
\widestnumber\key{MCS}

\ref\key Ad
\by C. Adams
\paper Volumes of $N$-cusped hyperbolic 3-manifolds
\jour J. London Math. Soc.
\vol 38 
\yr 1988
\pages 555-565
\endref

\ref\key AM
\by S. Akbulut and J. McCarthy
\book Casson's invariant for oriented homology $3$-spheres
 Mathematical Notes, 36.
\bookinfo Princeton University Press
\yr 1990
\endref

\ref\key Br
\by M. Brittenham
\paper Bounding canonical genus bounds volume
\paperinfo preprint
\endref

\ref\key GT
\by H. Goda and M. Teragaito
\paper Tunnel number one genus one non-simple knots
\paperinfo to appeat in Tokyo J. Math.
\endref

\ref\key HT
\by A. Hatcher and W. Thurston
\paper Incompressible surfaces in 2-bridge knot complements
\jour Inventiones Math. 
\vol 79 
\yr 1985 
\pages 225-246
\endref

\ref\key Ki
\by R. Kirby
\paper Problems in low-dimensional topology
\inbook Geometric topology (Athens, GA, 1993), Studies in Advanced Mathematics
\vol 2
\publ Amer. Math. Soc.
\publaddr Providence, RI
\yr 1997
\pages 35-473
\endref
\ref\key KK
\by M. Kobayashi and T. Kobayashi
\paper On canonical genus and free genus of a knot
\jour J. Knot Thy. Ram.
\vol 5
\yr 1996
\pages 77-85
\endref

\ref\key MCS
\by W. Magnus, A. Karrass,  and D. Solitar
\book Combinatorial group theory
\bookinfo John Wiley and Sons, Inc.
\yr 1966
\endref

\ref\key NZ
\by W. Neumann and D. Zagier
\paper Volumes of hyperbolic three-manifolds
\jour Topology 
\vol 24 
\yr 1985 
\pages 307-332
\endref

\ref\key Oz
\by M. Ozawa
\paper Satellite knots of free genus one
\paperinfo preprint
\endref

\ref\key Ro
\by D. Rolfsen 
\book Knots and Links 
\bookinfo Publish or Perish Press
\yr 1976
\endref

\ref\key Se
\by H. Seifert
\paper \"Uber das Geschlecht von Knoten 
\jour Math Annalen
\vol 110
\yr 1934
\pages 571-592
\endref


\ref 
\key Th1
\by  W. Thurston 
\book The Geometry and Topology of 3-manifolds
\bookinfo notes from lectures at Princeton University, 1978-80
\endref

\ref\key Th2
\bysame
\paper Three-dimensional manifolds, Kleinian groups and hyperbolic geometry 
\jour Bull. Amer. Math. Soc.
\vol 6 
\yr 1982
\pages 357-381
\endref

\ref\key We
\by J. Weeks
\paper SnapPea, a program for creating and studying hyperbolic 3-manifolds
\paperinfo available for download from www.geom.umn.edu
\endref

\endRefs
\enddocument